\documentclass{ws-p8-50x6-00}
\usepackage{amsmath,amscd,amssymb,amsthm}

\newtheorem{theorem}{Theorem}[section]
\newtheorem{corollary}[theorem]{Corollary} 
\newtheorem{lemma}[theorem]{Lemma}

\newtheorem{proposition}[theorem]{Proposition}
\theoremstyle{definition}
\newtheorem{definition}[theorem]{Definition}

\theoremstyle{remark}
\newtheorem{remark}[theorem]{Remark}

\newtheorem{example}[theorem]{Example}

\newcommand{\abs}[1]{\lvert#1\rvert}
\def\norm#1{\left\Vert#1\right\Vert}

\def\Q {{\mathbb Q}}
\def\I {{\mathbb I}}
\def\Pr{{\mathbb P}}
\def\C {{\mathbb C}}
\def\N{{\mathbb N}}
\def\e{{\varepsilon}}
\def\Z {{\mathbb Z}}
\def\I {{\mathbb{I}}}

\def\Un{{\mathcal U}\,}
\def\Nb{{\mathcal N}}
\def\Aut{{\mathrm{Aut}}\,}
\def\St{{\mathrm{St}}}

\def\UCB{{\mathrm{UCB}}\,}

\def\Homeo{{\mathrm{Homeo}}\,}
\def\RUCB{{\mathrm{RUCB}}\,}

\def\CB{{\mathrm{CB}}\,}
\def\GL{{\mathrm{GL}}\,}

\def\LO{{\mathrm{LO}}}
\def\SL{{\mathrm{SL}}\,}

\def\H{{\mathcal{H}}}

\def\1{{\mathbf 1}}
\def\Sa{{\mathcal S}}

\def\s{{\mathbb S}}

\begin{document}

\title{Remarks on actions on compacta 
by some infinite-dimensional groups}

\author{Vladimir Pestov}

\address{School of Mathematical and Computing Sciences,
Victoria University of Wellington, P.O. Box 600,
Wellington, New Zealand
\footnote{New permanent address beginning July 1, 2002: 
Department of Mathematics and Statistics,
University of Ottawa, Ottawa, Ontario, K1N 6N5, Canada.} \\
E-mail: vova@mcs.vuw.ac.nz \\
http://www.mcs.vuw.ac.nz/$^\sim$vova}


\maketitle

\abstracts{We discuss some techniques related to 
equivariant compactifications of uniform spaces and amenability
of topological groups. In particular, we give
a new proof of a recent result
by Glasner and Weiss describing the universal
minimal flow of the infinite symmetric group ${\mathfrak S}_\infty$
with the standard Polish topology, and extend
Bekka's concept of an amenable representation,
enabling one to
deduce non-amenability of the Banach--Lie groups 
$\GL(L_p)$ and $\GL(\ell_p)$, $1\leq p <\infty$.
}

\section{Introduction} 
Let a topological group act continuously by uniform isomorphisms
on a uniform space $X$. (One important situation
is where $X=G/H$ is a homogeneous
factor-space of $G$, equipped with the right
uniform structure.)
A compact space $K$, equipped with a continuous action of $G$,
is called an {\it equivariant compactification} of $G$ if there is
a uniformly continuous mapping $i\colon X\to K$ with dense image,
commuting with the action of $G$. 
Compactifications of this type always exist, moreover every
such $X$ admits a maximal $G$-equivariant compactification.

Here we discuss some ways in which 
equivariant compactifications can be used to
study minimal actions and amenability of some infinite-dimensional
groups. The latter term is used in an intuitive sense, to refer to
concrete topological groups of importance in mathematics, such as, for
instance, the full
unitary groups of the infinite-dimensional Hilbert spaces. 
Some of these groups form infinite-dimensional Lie groups 
in one or other sense.

A topological group $G$ is called amenable if every compact $G$-space
admits an invariant (regular Borel)
probability measure. In particular, $G$ is
{\it extremely amenable} if every compact $G$-space contains a
fixed point (that is, admits an invariant Dirac measure). 
No non-trivial locally compact group is extremely
amenable,\cite{V} but among infinite-dimensional groups extreme
amenability is not uncommon.\cite{Gr-M,Gl,P1,P4,G-P,P5}

A continuous action
of $G$ on a compact space $X$ is called minimal\cite{Aus} if 
the $G$-orbit of every point $x\in X$ is
everywhere dense in $X$. 
Every topological group $G$ possesses the 
universal minimal flow ($G$-space), ${\mathcal M}(G)$, such that
every other minimal $G$-flow is a factor of ${\mathcal M}(G)$. 
For locally compact groups the size of
the universal minimal flow is so immense that no
constructive description is ever likely. (Cf. e.g. \cite{Ell})
It comes as a surprise then that the universal minimal flow of at
least some infinite-dimensional groups is manageable. 

Moreover, it turns out that extremely amenable groups can be used
as a tool in order to give an explicit description of the universal
minimal flow ${\mathcal M}(G)$ even in cases where the flow is nontrivial.
If a topological group $G$ contains a `large' extremely amenable subgroup
$H$, then the universal minimal flow of $G$ is a subflow of the
equivariant compactification of the homogeneous space $G/H$, which is a much
smaller object than $G$ itself. In some cases, it enables one 
to describe ${\mathcal M}(G)$.
Such a technique was first used by the present author\cite{P1} 
in order to prove
that the circle $\s^1$ forms the universal minimal flow for the group of
orientation-preserving homeomorphisms of $\s^1$. 
Here we will use the
argument in order to give a more transparent proof of the
recent remarkable result by Glasner and Weiss,\cite{Gl-W}
who have characterized the universal
minimal flow of the infinite symmetric group 
${\mathfrak S}_\infty$, equipped with the
standard Polish topology, as the compact space of all linear orders on
$\N$. (The proof proposed here has an advantage that it extends
the result beyond the separable case, to groups of permutations of an
arbitrary infinite rank.)

Let us get back to the concept of an amenable topological group. 
A finer scale of `shades of amenability' is given by the following
concept:
say that a homogeneous factor-space
$G/H$ (or just a uniform $G$-space $X$)
is amenable in the sense of Eymard\cite{Ey} and Greenleaf\cite{Gre} if 
the maximal equivariant compactification of $G/H$ supports an
invariant probability measure. 

Here is an important particular case.
A unitary representation $\pi$ of a group $G$ in a Hilbert space
$\H$ is amenable in the sense of Bekka\cite{bekka}
if there is a state on the von Neumann algebra
of all bounded operators on $\H$, which is invariant under the
action of $G$ by conjugations. 
It turns out\cite{Pgafa} that a representation $\pi$ is amenable
if and only if the unit sphere in the Hilbert space
of the representation, upon which
$G$ acts by isometries, is an amenable uniform $G$-space.

In general, it is more difficult to verify
amenability of infinite-dimensional groups than that of
locally compact or discrete ones, because some tools
present in the locally compact case are missing.
For example, if a locally compact group $G$ contains a closed copy
of the free non-abelian group on two generators, then $G$ is
non-amenable, because amenability is inherited by closed subgroups
of locally compact groups. Not so beyond the locally compact case:\cite{dlH}
in fact, every topological group embeds into an extremely amenable
group as a topological subgroup.\cite{P4} Another example:
a locally compact group $G$ is amenable if and only if every
strongly continuous unitary representation of $G$ is amenable.\cite{bekka}
For infinite-dimensional groups, neither implication need hold.

Here we show that in some situations the property of amenability
is, in a sense, `partly' inherited by topological subgroups. 

We extend Bekka's concept as follows.
Say that a representation
$\pi$ of a group in a Banach space $E$
by bounded linear operators is amenable if
the projective space of $E$ (upon which the
group $G$ acts by isometries in a natural way) 
is an amenable $G$-space.

We show that every uniformly continuous
representation of an amenable topological group is amenable. Since
Eymard--Greenleaf amenability of an action (in particular, the Bekka
amenability of a representation) of a group $G$ is clearly inherited by
every subgroup $H<G$, we obtain a new possible way to prove that
a topological group $G$ is non-amenable: to find a uniformly
continuous representation $\pi$ of $G$ and a subgroup $H<G$ such that
the restriction of $\pi$ to $H$ is apriori non-amenable. 

The most natural class of infinite-dimensional groups admitting
uniformly continuous representations are Banach--Lie groups and
algebras of operators.
As an illustration of our methods, we show that
the general linear groups
$\GL(L_p)$ and $\GL(\ell_p)$, where $1\leq p<\infty$, 
are non-amenable if equipped with
the uniform operator topology. Even for Hilbert spaces
this seems to be a new result, answering a question
that Pierre de la Harpe asked me back in 1999.

\section{Some abstract nonsense}
\subsection{Uniformities and compactifications}
For a topological group $G$, we denote by $\Un_r(G)$ the
{\it Bourbaki-right} (= {\it Ellis-left}) uniform structure,
whose entourage basis consists of the sets
\[V_r=\{(x,y)\in G\times G\mid xy^{-1}\in V\},\]
and $V$ runs over the neighbourhood filter, $\Nb_G$, of $G$ at 
the neutral element $e_G$. 
The symbol $\RUCB(G)$ will denote the $C^\ast$-algebra of all
Bourbaki right uniformly continuous bounded complex-valued functions
on $G$, equipped with the supremum norm. 

Denote by $\Sa(G)$ the Samuel compactification of the uniform
space $(G,\Un_r(G))$, that is, the maximal ideal space of $\RUCB(G)$. 
This object
(together with the distinguished point, $e=e_G$)
is the well-known greatest ambit of $G$. (\cite{Aus,Br,P2})
In other words, $\Sa(G)$ is a $G$-ambit (a compact $G$-space with
a distinguished point having dense orbit), admitting a continuous
equivariant map, preserving the distinguished points, 
to any other $G$-ambit. 

Any two minimal compact $G$-subspaces of $\Sa(G)$ (whose existence
is guaranteed by Zorn's lemma) are isomorphic as $G$-spaces.\cite{Aus}
(This is a non-trivial fact, because there is, in general,
no {\it canonical} isomorphism.) This unique minimal $G$-space 
is denoted ${\mathcal M}(G)$ and called the {\it universal minimal
$G$-space} (or {\it $G$-flow}).

Let $H$ be a (closed or not) subgroup of a topological group $G$.
The {\it Bourbaki-right} uniform structure $\Un_r(G/H)$
is by definition the finest uniform structure on $G/H$ making the
factor-map 
\[G\ni g\mapsto gH\in G/H\]
uniformly continuous if $G$ is equipped with the uniformity $\Un_r(G)$.
In general, the uniformity $\Un_r(G/H)$ need not be separated even if
$H$ is a closed subgroup, and the topology generated on $G/H$ by
$\Un_r(G/H)$ may be coarser than the factor-topology on $G/H$.

The standard action of $G$ on $G/H$ on the left extends to the action
of $G$ on the Samuel compactification $\sigma(G/H,\Un_r(G/H))$.
(Notice that the Samuel compactification is
always a separated uniform space, and so the compactification map
need not be an embedding.) We will denote the latter compact space
by $\Sa_H(G)$.

The Banach $G$-module $C(\Sa_H(G))\cong \UCB(G/H,\Un_r(G/H))$
embeds into the Banach $G$-module $\UCB(\Sa(G))$.
Since the action of $G$ on the latter is well-known to be continuous,
the same is true of the action of $G$ on the former Banach space (and
$C^\ast$-algebra), and as a corollary, the action 
of $G$ on the compact space $\Sa_H(G)$ is continuous. 
With the image of the coset $H$ as the distinguished point,
$\Sa_H(G)$ is thus a $G$-ambit.

\subsection{Amenable groups and homogeneous spaces}

A topological group $G$ is called {\it amenable} 
if one of the
following equivalent conditions holds. All measures are assumed to be
regular Borel.
\begin{enumerate}
\item\label{one}
There is a
left-invariant mean on the space $\RUCB(G)$. 
\item \label{two}
There is an invariant probability measure on the greatest
ambit $\Sa(G)$. 
\item \label{three}
There is an invariant probability measure on every
compact space upon which $G$ acts continuously. 
\item There is an invariant probability measure on 
${\mathcal M}(G)$.
\end{enumerate}

See e.g. (\cite{Aus}, Chapter 12).
\smallskip

A topological group $G$ is called {\it extremely amenable} if 
one of the following equivalent conditions is true.

\begin{enumerate}
\item There is a multiplicative left-invariant mean on
$\RUCB(G)$.
\item There is a fixed point in $\Sa(G)$.
\item There is a fixed point in every compact space upon which
$G$ acts continuously. (The fixed point on compacta property.)
\item The universal minimal flow ${\mathcal M}(G)$ is
a singleton.
\end{enumerate}
Even if this property looks
exceedingly strong (in particular, no non-trivial locally compact group can
possess it\cite{V}), now we know numerous examples and entire classes of
infinite-dimensional groups that are extremely amenable. 
The following list is not exhaustive: 
the unitary group of an infinite-dimensional Hilbert
space with the strong operator topology;\cite{Gr-M} 
the group of classes of measurable maps from the unit interval to
the circle rotation group,\cite{F-W,Gl}
or, more generally, to any amenable locally compact group,\cite{P4}
equipped with the topology of convergence in measure;
the group of homeomorphisms of the closed (or open) unit interval with the
compact-open topology;\cite{P1}
the group of measure-preserving transformation of
the standard Lebesgue measure space with the weak topology, as well
as the group of measure class preserving transformations;\cite{G-P} the
group of isometries of the Urysohn universal metric space;
\cite{P4} 
unitary groups of certain von Neumann algebras and 
$C^\ast$-algebras.\cite{G-P} 
\smallskip

If $H$ is a subgroup of a topological group $G$, then the
homogeneous space $G/H$ (or the pair $(G,H)$) is
{\it Eymard--Greenleaf amenable}\cite{Ey,Gre} 
if there is a left-invariant mean
on the space $\UCB(G/H,\Un_r(G/H))$. Equivalently, there exists an
invariant probability measure on the ambit $\Sa_H(G)$.

More generally, one can talk of amenability of an action of
a group $G$ on a uniform space $X$ by uniform
isomorphisms. In such a situation, the topology on $G$ becomes
irrelevant.
\medskip

\begin{definition} Let a group $G$ act by uniform isomorphisms on a
uniform space $X$. Say that the action of $G$ is 
{\it Eymard--Greenleaf amenable,}
or that $X$ is an {\it Eymard--Greenleaf 
amenable uniform $G$-space}, if there 
exists a $G$-invariant mean on the space $\UCB(X)$. Equivalently
(by the Riesz representation theorem),
there exists an invariant probability measure on the Samuel
compactification $\sigma X$.
\label{def-amen}
\end{definition}
\smallskip

For example, in the case $X=(G/H,\Un_r(G/H))$ the above notion coincides
with Eymard--Greenleaf amenability. 

The following simple observation lends the concept 
some gravitas.
\smallskip

\begin{proposition} 
Every continuous action of an amenable locally compact
group $G$ on a uniform space $X$ by uniform isomorphisms is amenable.
\label{amenamen}
\end{proposition}

\begin{proof} Choose a point $x_0\in X$ and set, for each
$f\in \UCB(X)$ and every $g\in X$, 
\[\tilde f(g):=f(gx_0).\]
The function $\tilde f\colon G\to \C$ so defined is bounded
(obvious) and continuous, as the composition of two
continuous maps: the orbit map $g\mapsto gx_0$ and the 
function $f\colon X\to\C$. Also, for each $h\in G$,
\begin{eqnarray*}
\widetilde{^hf}(g) &=& ^hf(gx_0) \\
&=& f(h^{-1}gx_0) \\
&=& \tilde f(h^{-1}g) \\
&=& \,\,^h\tilde f(g),
\end{eqnarray*}
that is, the operator 
\[\alpha\colon \UCB(X)\ni f\mapsto\tilde f\in \CB(G)\]
is $G$-equivariant. (Here $\CB(G)$ denotes the $C^\ast$-algebra of
all bounded complex-valued continuous functions on $G$.)
It is also clear that $\alpha$ is 
positive, linear, bounded of norm one, and
sends the function $\mathbf 1$ to $\mathbf 1$.
Since $G$ is amenable and locally compact, there exists a left-invariant
mean $\phi$ on the space $\CB(G)$, and the composition
$\phi\circ\alpha$ is a $G$-invariant mean on $\UCB(X)$.
\end{proof}

This result is no longer true for more general topological groups,
cf. a discussion in Subsection \ref{reps}.

\subsection{More on the ambit $\Sa_H(G)$}

Let $G$ act continuously on a compact space $X$. Suppose there is a
point $\xi\in X$ stabilized by $H$. The orbit map 
\[G\ni g\mapsto g\xi\in X\]
then factors through the factor-space $G/H$, because for each
$h\in H$ one has $(gh)\xi= g(h\xi)=g\xi$. Denote the resulting
map $G/H\to X$ by $i$.
Since the orbit map
$G\to X$ is uniformly continuous relative to the uniformity 
$\Un_r(G)$, the inductive definition of the uniformity $\Un_r(G/H)$ implies
that $i$ is uniformly continuous as well. 
Consequently, $i$ extends in a unique way to a continuous 
equivariant map $\Sa_H(G)\to X$.
We conclude that $\Sa_H(G)$, with the distinguished point
$H$ (the coset of $e_G$), is the universal compact $G$-ambit with the
property that $H$ stabilizes the distinguished point. 

In general, the compact $G$-space $\Sa_H(G)$ need not be minimal.
The corresponding examples are easy to construct. 
\smallskip

However, notice the following.
\medskip

\begin{lemma}
Let $G$ be a topological group, and let $H$ be a closed subgroup.
Suppose the topological group $H$ is extremely amenable.
Then any minimal compact $G$-subspace, $\mathcal M$, of
$\Sa_H(G)$ is a universal minimal compact $G$-space.
\label{nonsense}
\end{lemma}

\begin{proof}
Let $X$ be an arbitrary minimal compact $G$-space.
Because of extreme amenability of $H$, there is a point
$\xi\in X$, stabilized by $H$. In view of the universality property
of $\Sa_H(G)$ described above, there is a morphism of $G$-spaces
$j\colon \Sa_H(G)\to X$ (taking $H$ to $\xi$). 
Because of minimality of $X$, the restriction of
the map $j$ to $\mathcal M$ is onto $X$. We are done.
\end{proof}

\begin{example} Let $G=\Homeo_+(\s^1)$, the group of orientation-preserving
homeomorphisms of the circle with the usual topology of uniform convergence,
and let $H$ be the isotropy subgroup of any chosen element
$\theta\in\s^1$. Then $H$ is isomorphic to the topological group
$\Homeo_+[0,1]$ and therefore extermely amenable.\cite{P1}
The right uniform 
factor-space $G/H$ is easily verified to be isomorphic to
the circle $\s^1$ with the unique compatible uniformity, and therefore
the ambit $\Sa_H(G)$ is $\s^1$ itself
with the distinguished point $\theta$. 
Since it is obviously a minimal $G$-space, we conclude by 
Lemma \ref{nonsense}: $\s^1$ is the universal minimal $\Homeo_+(\s^1)$-space.
This fact, established by the present author 
in (\cite{P1}), was probably the first instance where a non-trivial
universal minimal flow of any topological group has been computed
explicitely.
\end{example}
\smallskip

Here is another 
consequence of Lemma \ref{nonsense}, showing that
the class of extremely amenable group is closed under extensions,
similarly to the class of amenable groups. This result was
established (through a direct proof) during author's discussion
with Thierry Giordano and Pierre de la Harpe in April 1999, and is,
thus, a joint result.
\smallskip

\begin{corollary} Let $H$ be a closed normal subgroup of a topological
group $G$. If topological
groups $H$ and $G/H$ are extremely amenable, then
so is $G$.
\label{threeauthors}
\end{corollary}

\begin{proof}
In this case, the ambit $\Sa_H(G)$ is just the greatest ambit
$\Sa(G/H)$ of the factor-group, and it contains a fixed point
since $G/H$ is extremely amenable. Now we conclude by 
Lemma \ref{nonsense}.
\end{proof}

\section{The universal minimal flow of the 
infinite symmetric group\label{gl-we}}
Here we use Lemma \ref{nonsense} in order to
reprove a result by Glasner and Weiss\cite{Gl-W} describing the
universal minimal flow of the infinite symmetric group.
An idea of this new proof was briefly sketched by us in
(\cite{P5}, Exercises 11 and 12), but 
appears here in any detail for the first time.

Let $X$ be an infinite set (countable or not), and
let $G={\mathfrak S}_X$ denote 
the full group of permutations of $X$,
equipped with the topology of simple convergence on
$X$ viewed as a discrete space. For countable $X$, this topology is
well known to be Polish (separable completely metrizable).

Denote by $\LO_X$ the set of all linear orders on $X$,
equipped with the (compact) 
topology induced from $\{0,1\}^{X\times X}$. (Here a linear
order $\prec$ is identified with the characteristic function of
the corresponding relation $\{(x,y)\in X\times X\colon x\prec y\}$.)

The group 
${\mathfrak S}_X$
acts on $\LO_X$ by double permutations:
\[
(x\,\,^\sigma\!\!\prec y)\Leftrightarrow (\sigma^{-1}x\prec \sigma^{-1}y)\]
for all $\prec\,\,\in \LO_X$, $\sigma\in {\mathfrak S}_X$, and $x,y\in X$.
This action is
continuous and minimal (an easy check). 

A linear order $\prec$ on $X$ is called {\it $\omega$-homogeneous}
if every finite subset $A\subset X$ can be mapped
onto any other subset $B\subset X$ of the same cardinality by an
order-preserving bijection (order automorphism) of $(X,\prec)$.
In particular, it follows that $\prec$ is a dense order without least
and greatest elements. (In the case where $X$ is countable, this
condition is equivalent to $\omega$-homogeneity.) 

Every infinite set $X$
supports an $\omega$-homogeneous linear order. 
(Here is one proof: $X$ can be given the structure of an
ordered field, because it has 
the same cardinality as $\Q(X)$, the 
purely transcendental field extension of $\Q$, and the field
$\Q(X)$ is well known to be linearly orderable. 
And every linearly ordered field is $\omega$-homogeneous
due to the existence of piecewise-linear monotone maps.)
Choose an arbitrary such order on $X$, say
$\prec$. 

Let $H=\Aut(\prec)$ be the subgroup of all
permutations preserving the linear order $\prec$.
The left factor-space $G/H={\mathfrak S}_X/\Aut(\prec)$ can be
identified with a certain collection of linear orders on $X$, namely
those obtained from $\prec$ by a permutation. Denote this collection
by $\LO_\prec$. 
Thus, $G/H\cong \LO_\prec$ embeds into $\LO_X$.

As every compact space, $\LO_X$ supports a unique compatible uniform
structure. It induces a totally bounded uniform structure on $\LO_\prec$. 
\smallskip

\begin{lemma}
The uniform
structure on $G/H\cong \LO_\prec$, 
induced from the compact space $\LO_X$, 
coincides with the right uniform structure
$\Un_r({\mathfrak S}_X/\Aut(\prec))$.
\label{finest}
\end{lemma}

\begin{proof} We want to show that the uniform
structure on $G/H$, induced from the compact space $\LO_X$, is 
the finest uniform structure making the quotient map
\[
{\mathfrak S}_X\to {\mathfrak S}_X/\Aut(\prec)\cong\LO_\prec
\] 
right uniformly continuous. The proof consists of two parts.
\smallskip

(1) The map $\sigma\mapsto \,\,^\sigma\!\!\prec$ is uniformly continuous.
\smallskip

Let $F=\{x_1,\dots,x_n\}\subset\omega$ be 
any finite subset, determining the following
standard basic entourage of the uniformity of $\LO_X$: 
\[
W_F:=\{(<_1,<_2)\in \LO_X\times\LO_X\colon <_1\vert_F=<_2\vert_F\}.
\]
Denote by $\St_F$ the common isotropy subgroup of all
$x_i\in F$, that is,
\[
\St_F:=\{\sigma\in {\mathfrak S}_X\mid \sigma(x_i)=x_i,~~i=1,2,\dots,n\}.
\]
This $\St_F$ is an open subgroup of ${\mathfrak S}_X$ and in particular a
(standard) open neighbourhood of the identity. As such, it determines
an element of the Bourbaki-right uniformity $\Un_r({\mathfrak S}_X)$:
\[
V_F:=\{(\sigma,\tau)\in {\mathfrak S}_X\times {\mathfrak S}_X\mid
\sigma\tau^{-1}\in\St_F\}.
\]
In other words, $(\sigma,\tau)\in V_F$ iff for all $i=1,2,\dots,n$
one has $\tau^{-1}x_i=\sigma^{-1}x_i$. 

If now $(\sigma,\tau)\in V_F$, then for every $i,j=1,2,\dots,n$ one has
\begin{eqnarray*}
x_i\,\,^\sigma\prec x_j &\Leftrightarrow& \sigma^{-1}x_i\prec
\sigma^{-1}x_j \\
 &\Leftrightarrow& \tau^{-1}x_i\prec \tau^{-1}x_j \\ 
&\Leftrightarrow&
x_i\,\,^\tau\prec x_j,
\end{eqnarray*}
meaning that the restrictions of the orders
$^\sigma\prec$ and $^\tau\prec$ to $F$ coincide, and thus
$(^\sigma\prec,\,\,^\tau\prec)\in W_F$. 
\smallskip

(2) The image of the entourage $V_F$ under
the (Cartesian square of) the map $\sigma\mapsto\,\,^\sigma\prec$
is {\it exactly} all of $W_F\cap (\LO_\prec\times\LO_\prec)$. 
\smallskip

Indeed, suppose $(<_1,<_2)\in W_F\cap (\LO_\prec\times\LO_\prec)$,
that is, $<_1$ and $<_2$ are linear orders on $X$, obtained from
$\prec$ by suitable permutations, and
whose restrictions to a finite subset $F$ coincide. 

Choose two permutations $\sigma,\tau\in {\mathfrak S}_X$ such that
$<_1=\,\,^\sigma\prec$ and $<_2=\,\,^\tau\prec$. 
For each $i=1,2,\dots,n$, one necessarily has 
\[\sigma^{-1}x_i =\tau^{-1}x_i\]
(if it were not so, then the orders $\,\,^\sigma\prec$ and 
$\,\,^\tau\prec$ would differ on $F$). Consequently,
$(\sigma,\tau)\in V_F$. 
\smallskip

Now we conclude that every uniform structure, $\Un$, on $\LO_\prec$
that makes the map 
\[({\mathfrak S}_X,\Un_r)
\ni\sigma\mapsto\,\,^\sigma\prec\,\in(\LO_\prec,\,\Un)\] 
uniformly continuous,
must be coarser than the restriction of the uniformity of $\LO_X$
to $\LO_\prec$. Indeed, for every element $W\in\Un$ there is, by the
assumed uniform continuity of the above map,
a finite $F\subset\omega$
with the image of $V_F$ contained in $W$, that is, with $W_F\subseteq W$.
This accomplishes the argument. 
\end{proof}

\begin{lemma}
The ambit $\Sa_{\Aut(\prec)}({\mathfrak S}_X)$ is isomorphic
to $\LO_X$, with the distinguished element $\prec$.
\end{lemma}

\begin{proof}
By Lemma \ref{finest}, 
$({\mathfrak S}_X/\Aut(\prec),\Un_r({\mathfrak S}_X/\Aut(\prec))$
embeds into $\LO_X$ as a uniform subspace and 
an ${\mathfrak S}_X$-subspace. 
Also, $\LO_\prec$ is everywhere dense in $\LO_X$.
As a consequence, the Samuel compactification of the precompact
uniform space 
$({\mathfrak S}_X/\Aut(\prec),\Un_r({\mathfrak S}_X/\Aut(\prec))$
is simply its completion, that is, $\LO_X$. 
\end{proof}

 An application of Lemma \ref{nonsense}
(bearing in mind that the topological group
$H=\Aut(\prec)$ is extremely amenable \cite{P1}) yields immediately:
\vskip .3cm
\begin{theorem}[Glasner and Weiss \cite{Gl-W}]
The compact space $\LO_X$ forms the universal 
minimal ${\mathfrak S}_X$-space.
\qed
\end{theorem}
\smallskip

\begin{remark} 
The original theorem by Glasner and Weiss was 
established in the case of countable
$X$. Our proof remains true for symmetric groups of arbitrary
infinite rank.
\end{remark}
\smallskip

\begin{remark} The group ${\mathfrak S}_X$ contains, as a dense
subgroup, the union of the directed family of permutation subgroups
of finite rank, and consequently it is amenable. As a result, there
is an invariant probability measure on the compact set
$\LO_X$. Glasner and Weiss\cite{Gl-W} have proved that
such a measure is unique, that is, the action by ${\mathfrak S}_X$
on ${\mathcal M}({\mathfrak S}_X)\cong\LO_X$ is uniquely ergodic.
\smallskip

Their argument can be made quite elementary
(no Ergodic Theorem!) as follows. 
Let $\mu$ be a ${\mathfrak S}_X$-invariant probability measure on $\LO_X$.
If $F\subset X$ is a finite subset, then every linear order, $<$, on $F$
determines a cylindrical subset
\[C_<:=\{\prec\,\in\LO_X\colon \prec\vert_F=<\}\subset\LO_X.\]
Every two sets of this form, corresponding to different orders
on $F$, are disjoint and can be taken to each
other by a suitable permutation. As there are $n!$ of such sets,
where $n=\abs F$,
the $\mu$-measure of each of them 
must equal $1/n!$. Consequently, the functional $\int d\mu$ is
uniquely defined on the characteristic functions of cylinder sets
$C_<$, which functions are continuous and separate points, because
sets $C_<$ are open and closed and
form a basis of open subsets of $\LO_X$. Now the Stone--Weierstrass
theorem implies uniqueness of $\int d\mu$ on all of $C(\LO_X)$.
\end{remark}
\smallskip

\begin{remark}
Every extremely amenable subgroup $H$ of ${\mathfrak S}_X$ is contained in
one of the subgroups of the form $\Aut(\prec)$.
(Indeed, $H$ must possess a fixed point in the space $LO_X$,
that is, preserve a linear order $\prec$ 
on $X$.)

At the same time,
not every subgroup of the form $\Aut(\prec)$ is extremely amenable.
For example, if the linear order $\prec$ is such that
for some cover of $X$ by
three disjoint convex subsets $A,B,C$ one has $A<B<C$, $A$ and $C$
are densely ordered, and $B$ has type $\Z$, then the group
$\Aut(\prec)$ is topologically isomorphic to the product of three
groups of order automorphisms, and since $\Aut(B)\cong\Z$ is not
extremely amenable, neither is $\Aut(\prec)$.

On the other hand,
a similar construction can be used to produce examples of groups of
type $\Aut(\prec)$ which are extremely amenable even if the order
$\prec$ is not dense (admits gaps).
\label{maximal}
\end{remark}
\smallskip

\begin{example} 
The {\it tame topology} on the group $U(\infty) = \cup_{i=1}^\infty U(n)$
is the topology of simple convergence on the sphere 
$\s(\infty)=\cup_{i=1}^\infty\s^n$ (the intersection of the unit
sphere of $\ell_2$ 
with the direct limit space $\C^\infty$),
viewed as discrete. Thus,
$U(\infty)$ receives the subgroup topology from ${\mathfrak S}_{\s(\infty)}$.
This topology is of
considerable interest in representation theory of the infinite
unitary group,\cite{Olsh} where unitary representations strongly
continuous with regard to the tame topology are called
{\it tame representations.}

As a consequence of the Remark \ref{maximal},
the group $U(\infty)$ with the tame topology is not extremely
amenable: indeed, it is easy to see that
no linear order on $\s(\infty)$ is preserved
by all operators from $U(\infty)$. Thus, the universal minimal flow
${\mathcal M}(U(\infty)_{tame})$ is nontrivial.
\label{tame}
\end{example}
\smallskip

\begin{remark}
Let a group $G$ act by uniform isomorphisms on a uniform space $X$.
The pair $(G,X)$ has the 
{\it Ramsey--Dvoretzky--Milman property} if for every
bounded uniformly continuous function $f$ from $X$ to a finite-dimensional
Euclidean space, every finite $F\subseteq X$, and each $\e>0$ 
there is a $g\in G$
such that the oscillation of $f$ on the translate $gF$ is less than $\e$.
This concept links extreme amenability with Ramsey theory, because
a topological group $G$ is extremely amenable if and only if every
continuous transitive
action of $G$ by isometries on a metric space has the 
Ramsey--Dvoretzky--Milman property.\cite{P5} 

The statement in Example \ref{tame} can be strengthened:
a result by Graham\cite{Grh1} on the so-called sphere-Ramsey spaces
implies that already
the pair $(U(\infty),\s(\infty))$, where $\s(\infty)$ is
equipped with the discrete ($\{0,1\}$-valued) metric, does not have
the Ramsey--Dvoretzky--Milman property. 

This sort of dynamical properties, formulated for appropriate groups
of affine transformations,  is 
linked to the central open
question of Euclidean Ramsey theory: is every finite spherical metric space
Ramsey? \cite{Grh} 
\end{remark}

\section{Amenable representations} 

\subsection{The projective space}

Let $E$ be a (complex or real) Banach space.
Denote by $\Pr_E$ the projective space of $E$.
If we think of $\Pr_E$ as a factor-space of the unit
sphere $\s_{E}$ of $E$, then $\Pr_E$ becomes a metric space via
the rule
\[d(x,y) = \inf\{\norm{\xi-\zeta}\colon \xi,\zeta\in \s_{\H},
p(\xi)=x, p(\zeta)=y\},\]
where $p\colon \s_{E}\to\Pr_E$ is the canonical factor-map. Notice
that the infimum in the formula above is in fact minimum. The proof of
the triangle inequality is based on the invariance of the norm distance
on the sphere under multiplication by scalars. The above metric on the
projective space is complete.
\smallskip

Let $T\in\GL(E)$ be a bounded linear invertible operator on
a Banach space $E$. Define a mapping $\tilde T$ from the projective
space $\Pr_E$ to itself as follows: for every $\xi\in \s_E$ set
\[\tilde T(p(\xi)) = p\left(\frac{T(\xi)}{\norm{T(\xi)}}\right).\]
The above definition is clearly
independent on the choice of a representative,
$\xi$, of an element of the projective space $x\in\Pr_E$. 
\smallskip

\begin{lemma} The mapping $\tilde T$ is a uniform isomorphism
(and even a bi-Lipschitz isomorphism)
of the projective space $\Pr_E$.
\label{bilip}
\end{lemma}

\begin{proof}
It is enough to show that $\tilde T$ is uniformly continuous,
because $\widetilde{TS}=\tilde T\tilde S$ and so
$\widetilde{T^{-1}}=\tilde T^{-1}$. Let
$x,y\in\Pr_E$, and let $\xi,\zeta\in\s_E$ be such that
$p(\xi)=x$, $p(\zeta)=y$, and $\norm{\xi-\zeta}=d(x,y)$. 
Both $\norm{T(\xi)}$ and  $\norm{T(\zeta)}$ are bounded below
by $\norm{T^{-1}}^{-1}$, and therefore
\begin{eqnarray*}
d(\tilde T(x),\tilde T(y)) &\leq & \frac{\pi}2\norm{T^{-1}}
\cdot\norm{T(\xi)-T(\zeta)}\\
&\leq & \frac{\pi}2\norm{T^{-1}}\cdot\norm T\cdot \norm{\xi-\zeta}
\\ &= &  \frac{\pi}2\norm{T^{-1}}\cdot \norm T d(x,y).
\end{eqnarray*}
\end{proof}

Let us recall the following notion from theory of transformation
groups.

\begin{definition} Let a group $G$ act by uniform isomorphisms on
a uniform space $X=(X,\Un_X)$. The action is called
{\it bounded} (or else
{\it motion equicontinuous}) if for every $U\in\Un_X$ there is
a neighbourhood of the identity, $V\ni e_G$, such that
$(x,g\cdot x)\in U$ for all $g\in V$ and $x\in X$. 
\end{definition}
\smallskip

Notice that every bounded action is continuous.
\smallskip

\begin{example} 
The action of $\GL(E)$ on the unit sphere $\s_E$ (and moreover the
unit ball) of a Banach
space $E$ is bounded, by the very definition of the uniform
operator topology.
\end{example} 
\smallskip

\begin{lemma} 
\label{correspo}
The correspondence 
\[\GL(E)\ni T\mapsto \tilde T\]
determines an action of the general linear group $\GL(E)$ on
the projective space $\Pr_E$
by uniform isomorphisms. With respect to
the uniform operator topology on $\GL(E)$, the action is bounded.
\end{lemma}

\begin{proof} The first part of the statement is easy to check
using Lemma \ref{bilip}.
As to the second, if $\norm{T-\I}<\e$,
then for every $\xi\in\s_E$ 
\begin{eqnarray*}
\norm{\tilde T(x)-x}&\leq& \frac{\pi}2\norm{T^{-1}}\cdot\norm{T(x)-x}
\\
&<& \frac{\pi\e}{2(1-\e)}.
\end{eqnarray*}
\end{proof}

\subsection{Extension of Bekka's amenability}
We want to reformulate the concept of an
amenable representation in the sense of Bekka in order to present
a natural extension of it.

Let $\pi$ be a unitary representation of a group $G$ in a
Hilbert space $\H$. One says that $\pi$ is {\it amenable} \cite{bekka} 
if there exists a state,
$\phi$, on the von Neumann algebra ${\mathcal B}(\H)$
of all bounded operators on the space $\H$ of representation, which is
invariant under the action of $G$ by inner automorphisms: 
$\phi(\pi(g)T\pi(g)^{-1})=\phi(T)$ for
every $T\in B(\H)$ and every $g\in G$. 

The group $G$ acts on the unit sphere $\s_{\H}$ by isometries, and 
it was shown by the author\cite{Pgafa} that a unitary representation
$\pi$ of a group $G$ in a Hilbert space $\H$ is amenable if and only
if $\s_{\H}$ is an amenable uniform $G$-space in the sense of our
Definition \ref{def-amen}, that is,
there exists a $G$-invariant mean on the space $\UCB(\s_{\H})$ or,
equivalently, an invariant probability measure on the Samuel
compactification of the sphere $\s_{\H}$.

While the implication $\Rightarrow$ is based on some results obtained
by Bekka using 
deep techniques by Connes, the implication
$\Leftarrow$ is elementary. We need to reproduce it here.
\smallskip

$\triangleleft$ Let $\psi$ be a $G$-invariant mean on $\UCB(\s_\H)$.
Every bounded linear operator $T$ on
$\H$ defines a bounded uniformly continuous (in fact, even
Lipschitz) function $f_T\colon\s_\H\to\C$ by the rule
\[
\s_\H\ni\xi\mapsto f_T(\xi):=\langle T\xi, \xi\rangle\in \C.
\]
Now set $\phi(T):=\psi(f_T)$. This $\phi$ is a $G$-invariant
state on ${\mathcal B}(\H)$. $\triangleright$
\smallskip

Notice that the function $f_T$ in the proof 
above is symmetric: 
for every $\lambda\in \C$, $\abs\lambda =1$, and each 
$\xi\in\s_\pi$, one has $f_T(\lambda \xi)=f_T(\xi)$. In other words,
$f_T$ is constant on the preimages of $p$. 
Consequently,
$f_T$ factors through a function $\tilde f_T$ on the projective space
$\Pr_{\H}$; clearly, $\tilde f_T$ is also uniformly continuous and 
bounded. It means that
the above proof only uses the existence of
a $G$-invariant mean on the function space $\UCB(\Pr_{\H})$. 

On the other hand, the Banach space (and $G$-module) $\UCB(\Pr_{\H})$
admits an obvious equivariant embedding into $\UCB(\s_{\H})$; namely,
it can be identified with the Banach $G$-submodule of all functions
symmetric in the above sense. The restriction of a $G$-invariant mean
from $\UCB(\s_{\H})$ to $\UCB(\Pr_{\H})$ is again a $G$-invariant mean.

We have thus established the following.
\smallskip

\begin{theorem} 
\label{criter-amen}
A unitary representation $\pi$ of a group $G$ in a
Hilbert space $\H$ is amenable if and only if the projective
space $\Pr_{\H}$ is an Eymard--Greenleaf amenable uniform $G$-space.
\qed
\end{theorem}
\smallskip

The advantage of this reformulation is that it allows for an
extension of the concept of an amenable representation to group
representations by bounded linear operators that are not necessarily
unitary. 
\smallskip

\begin{definition}
Say that a representation $\pi$ of a group $G$ by bounded linear
operators in a normed space $E$ is {\it amenable} if the action of $G$
by uniform isometries on the projective space $\Pr_E$, associated
to $\pi$ as in Lemma \ref{correspo}, 
is Eymard--Greenleaf amenable in the sense of 
Definition \ref{def-amen}.
\end{definition}
\smallskip

\begin{theorem} Let $G$ be a locally compact group, and let
$\mu$ be a quasi-invariant measure on $G$. Let $1\leq p<\infty$. 
The left quasi-regular representation of $G$ in $L_p(\mu)$ is
amenable if and only if $G$ is amenable. 
\label{quasi}
\end{theorem}

\begin{proof}
 The left quasi-regular representation, $\gamma$,
of $G$ in $L_p(\mu)$ (cf. e.g. \cite{Wa}) is given by the formula
\[^gf(x) = \left(\frac{d(\tau\circ g^{-1})}{d\tau}\right)^{\frac 1p}
f(g^{-1}x),\]
where $d/d\tau$ is the Radon-Nykodim derivative. It is a strongly
continuous representation by isometries.
Necessity
($\Rightarrow$) thus follows at once from Proposition \ref{amenamen}.

To prove sufficiency ($\Leftarrow$), assume $\gamma$ is
amenable. Then there exists an invariant mean, $\phi$, on $\UCB(\s_p)$,
where $\s_p$ stands for the unit sphere in $L_p(\mu)$. 
For every Borel subset $A\subseteq G$, 
define a function $f_A\colon\s_p\to\C$
by letting for each $\xi\in\s_p$
\[f_A(\xi) = \norm{\xi\cdot\chi_A}^p, \]
where $\chi_A$ is the characteristic function of $A$. 
The function $f_A$ is bounded and uniformly continuous on $\s_p$.
For every $g\in G$, 
\begin{eqnarray*}
^gf_A(\xi) &=& f_A\left(^{g^{-1}}\xi\right) \\[2mm]
&=& \int_A\left\vert ^{g^{-1}}\xi(x)\right\vert^p d\mu(x) \\[2mm]
&=& \int_A \frac{d\mu\circ g}{d\mu}\left\vert\xi(gx)
  \right\vert^p d\mu(x) \\[2mm]
&=&
\int_{gA}\left\vert\xi(y)\right\vert^p d\mu(y)
\\[2mm]
&=& f_{gA}(\xi),
\end{eqnarray*}
that is, $^gf_A=f_{gA}$. 
It is now easily seen that $m(A):=\phi(f_A)$ is a finitely additive 
left-invariant measure on $G$, vanishing on locally null sets, and so
$G$ is amenable. 
\end{proof}

\subsection{\label{reps}Uniformly continuous
representations}
In contrast with Proposition \ref{amenamen},
even a strongly continuous unitary representation of an amenable
non-locally compact topological group need not be amenable.
The simplest example is the standard representation of the full
unitary group $U(\H)_s$ of an infinite-dimensional Hilbert
space, equipped with the strong topology.
It is not amenable because it contains, as a subrepresentation, the
left regular representation of a free nonabelian group, which is
not amenable. (The left regular representation of a locally
compact group $G$ is amenable if and only if $G$ is 
amenable,\cite{bekka} cf. also Theorem \ref{quasi}.)

A part of the story here
is that when a topological group $G$ acts
continuously by uniform isomorphisms on a uniform space $X$,
the resulting representation of $G$ by
isometries in $\UCB(X)$ need not be continuous. 
(This is the case,
for instance, in the same example $G=U(\H)_s$, $X=\s_{\H}$.) 
Equivalently,
the extension of the action of $G$ to the
Samuel (uniform) compactification $\sigma X$ is discontinuous, and
therefore one cannot deduce the existence of an invariant measure
on $\sigma X$ from the assumed amenability of $G$.

Here is a simple case where the continuity of action is assured.
\smallskip

\begin{lemma} Suppose a topological group $G$ acts in a bounded way
on a uniform space $X$. Then the resulting representation of $G$
in $\UCB(X)$, as well as the action of $G$ on $\sigma X$,
are both continuous.
\end{lemma}

\begin{proof} 
Since $G$ acts on $\UCB(X)$ by isometries, it is enough to show
that the mapping $G\ni g\mapsto \,^gf\in\UCB(X)$ is continuous
at identity. By a given $\e>0$, choose a $U\in\Un_X$ using the
uniform continuity of $f$, and a symmetric neighbourhood 
$V\ni e_G$ so that
$\abs{f(x)-f(y)}<\e$ whenever $(x,y)\in U$ and
$(x,g\cdot x)\in U$ once $g\in V$ and $x\in X$. 
Now for each $x\in X$, 
$\abs{f(x)-f(g^{-1}x)} \leq \e$ once $g\in V$, that is,
$\norm{f-\,^gf}_{sup}<\e$, and we are done. 

Now recall that the Samuel (maximal uniform) compactification of
$X$ is the maximal ideal space of $\UCB(X)$. It is a simple observation
(which was made, for instance, by Teleman\cite{Tel}) 
that a representation of
a topological group by isomorphisms of a commutative 
$C^\ast$-algebra is strongly
continuous if and only if the associated action of $G$ on the maximal
ideal space is continuous. 
\end{proof} 

\smallskip
The following three corollaries are immediate.
\smallskip

\begin{corollary}
Let a topological group $G$ act in a bounded way on a uniform space $X$.
If $G$ is amenable, then $X$ is an Eymard--Greenleaf amenable uniform
$G$-space. \qed
\end{corollary}
\smallskip

\begin{corollary} 
Let $\pi$ be a uniformly continuous
representation of a topological group in a Banach space $E$.
If $G$ is amenable, then $\pi$ is an amenable representation.
\qed
\end{corollary}
\smallskip

\begin{corollary} 
Let $E$ be a Banach space, and let $G$ be a topological subgroup
of $\GL(E)$ (equipped with the uniform operator topology).
If $H$ is a subgroup of $G$ and the restriction of the standard
representation of $\GL(E)$ in $E$
to $H$ is non-amenable, then $G$ is a non-amenable
topological group. \qed
\label{tretie}
\end{corollary}

\subsection{Groups of operators}

\begin{theorem} The general linear groups $\GL(L_p)$ and
$\GL(\ell_p)$,  $1\leq p<\infty$, 
with the uniform operator topology are non-amenable.
\end{theorem}

\begin{proof} 
Both spaces $L_p$ and $\ell_p$ can be realized as $L_p(\mu)$,
where $\mu$ is a quasi-invariant measure on a non-amenable
locally compact group, $H$. (For instance, $H=\SL(2,\C)$ for 
the continuous case and $H=\SL(2,\Z)$ for the purely atomic one.) 
Identify $H$ with an (abstract, non-topological) 
subgroup of $\GL(L_p(\mu))$ via the left quasi-regular representation,
$\gamma$. 
The restriction of the standard representation of $\GL(L_p(\mu))$
to $H$ is $\gamma$, which is a non-amenable representation
(Theorem \ref{quasi}), and Corollary \ref{tretie} applies.
\end{proof} 

It is interesting to compare the above result with the following.
\smallskip

\begin{theorem} 
The isometry group
${\mathrm{Iso}}(\ell_p)$, $1\leq p<\infty$,
$p\neq 2$, equipped with the strong operator
topology, is amenable, but not extremely amenable.
\label{l}
\end{theorem}

\begin{proof} The isometry groups in question, as abstract groups,
have been described by Banach in his classical 1932 treatise
\cite{banach} 
(Chap. XI, \S 5, pp. 178--179).
For $p>1$, $p\neq 2$ the
group ${\mathrm{Iso}}(\ell_p)$ is isomorphic to the semidirect
product of the group of permutations ${\mathfrak S}_\infty$
and the countable power $U(1)^\N$ (in the complex case) or
$\{1,-1\}^\Z$ (in the real case). Here the group of permutations
acts on $\ell_p$ by permuting coordinates, while the group 
of sequences of scalars of absolute value one acts by coordinate-wise
multiplication. The semidirect product is formed with regard to an
obvious action of ${\mathfrak S}_\infty$ on 
$U(1)^{\N}$ (in the real case, $\{1,-1\}^{\N}$).

The strong operator topology restricted to the group 
${\mathfrak S}_\infty$ is the standard Polish topology,
and restricted to the product
group, it is the standard product topology. Thus, 
${\mathrm{Iso}}(\ell_p)\cong {\mathfrak S}_\infty\ltimes U(1)^\N$
(correspondingly, ${\mathfrak S}_\infty\ltimes \{1,-1\}^{\N}$)
is the semidirect product of a Polish group with a compact metric
group. Since ${\mathfrak S}_\infty$ is an amenable topological
group, so is ${\mathrm{Iso}}(\ell_p)$. Since the non-extremely
amenable group ${\mathfrak S}_\infty$ is a topological 
factor-group of ${\mathrm{Iso}}(\ell_p)$, the latter group
is not extremely amenable either.
\end{proof}

\begin{remark}
Using the description of the universal minimal flow of
${\mathfrak S}_\infty$ due to Glasner and Weiss discussed in 
Section \ref{gl-we}, as well as some
standard means of uniform topology,
one can show that the universal minimal flow of the 
topological group ${\mathrm{Iso}}(\ell_p)$ is homeomorphic
to the product of the compact space $\LO_{\infty}$ of all
linear orders on the natural numbers with the compact
group $U(1)^{\N}$ (complex case) or $\{1,-1\}^{\N}$ (real case).
This compact space is equipped with a skew product
action: 
\[(\sigma,f)\cdot (\prec,g) = (^{\sigma}\prec,\, f\cdot\,^{\sigma}g),
\]
where $^{\sigma}g(n) = g(\sigma^{-1}n)$. Again, this action is
uniquely ergodic. We leave the details out.
\end{remark}
\smallskip

\begin{remark}
By contrast, for $p=2$ the unitary group of an infinite-dimensional
Hilbert space with the strong operator topology is extremely
amenable. This is due to Gromov and Milman\cite{Gr-M}.

We conjecture that the groups ${\mathrm{Iso}}(L_p)$,
$1\leq p<\infty$, with
the strong operator topology are all extremely amenable.
\end{remark}
\smallskip

\section*{Acknowledgments}
I am grateful to Joshua Leslie and Thierry Robart for their
hospitality during the conference
on Infinite Dimensional Lie Groups in Geometry and Representation
Theory and for their patience during the preparation of this volume.
Stimulating conversations with Thierry Giordano,
Eli Glasner, Pierre de la
Harpe, and Michael Megrelishvili are acknowledged.
The present research has been supported by the Marsden Fund of the Royal
Society of New Zealand through the grant project
`Geometry of high-dimensional structures: dynamical aspects.'


\end{document}